\documentclass[12pt]{article}
\usepackage[latin1]{inputenc}
\usepackage{amssymb,amsmath,amsfonts}
\newtheorem{theorem}{Theorem}
\newtheorem{proposition}[theorem]{Proposition}
\newtheorem{lemma}[theorem]{Lemma}

\newtheorem{corollary}[theorem]{Corollary}
\newtheorem{fact}[theorem]{Fact}

\newcommand{\R}{\mathbb{R}}
\newcommand{\Q}{\mathbb{Q}}
\newcommand{\Sf}{\mathbb{S}}
\newcommand{\C}{\mathbb{C}}

\newcommand{\spa}{\mbox{span}}

\newcommand{\CP}{\mathbb{CP}}
\newcommand{\po}{{\hspace*{-1ex}}{\bf .  }}

\newcommand{\nap}{\nabla^{\perp}}

\newcommand{\E}{{\cal E}}
\def\<{{\langle}}
\def\>{{\rangle}}

\def\a{\alpha}

\def\bea{\begin{eqnarray*} }
\def\eea{\end{eqnarray*} }
\def\be{\begin{equation} }
\def\ee{\end{equation} }

\def\nap{\nabla^\perp}
\def\proof{\noindent\emph{Proof: }}
\def\qed{\ifhmode\unskip\nobreak\fi\ifmmode\ifinner
\else\hskip5 pt \fi\fi\hbox{\hskip5 pt \vrule width4 pt
height6 pt  depth1.5 pt \hskip 1pt }}

\begin{document}

\title{Isometric deformations of  isotropic  surfaces}
\author{M.\ Dajczer and Th.\ Vlachos}
\date{}
\maketitle
%\subjclass[2010]{Primary 53C40; Secondary 53B25.}
%\keywords{}

\begin{abstract} It was shown by Ramanathan \cite{R} that any compact 
oriented minimal surface in the three-dimensional round 
sphere admits at most a finite set of pairwise noncongruent minimal 
isometric immersions. In this paper, we extend this result
to isotropic surfaces in spheres of arbitrary dimension.  
The case of non-compact isotropic surfaces in space forms is also 
addressed. 
\end{abstract}

It is standard fact that a simply-connected minimal surface 
$g\colon L^2\to\Q^3$  lying in  a complete  simply-connected three-dimensional 
space form possess an associated  one-parameter family of minimal isometric 
deformations $g_\theta\colon L^2\to\Q^3$, $\theta\in\Sf^1=[0,\pi)$.
On the other hand, if $L^2$ is non-simply-connected then the set of all 
minimal isometric deformations turns out to be either finite or still the 
full circle $\Sf^1$.  

Of particular interest is when the ambient space is a sphere since only 
in that case the minimal surface $g\colon L^2\to\Sf^3$ is allowed to be compact. 
It was proved by Ramanathan \cite{R} that a compact  surface only allows 
finitely many minimal isometric deformations. 
Moreover, the surface must be equivariant, that is, any intrinsic 
one-parameter family of isometries becomes extrinsic, as is the case
of the Clifford torus.

Our goal in this paper is to extend Ramanathan's results to any codimension. 
But, as already  shown by him this can only be achieved under some 
strong additional condition on the surface.  From the known results listed 
in \S $2$, we can see why being isotropic comes in naturally for that purpose. 
In order to  explain what for a surface being isotropic means, 
we need the  concepts of higher fundamental forms and curvature ellipses 
that we briefly recall next and refer to \S $1$ for details. 

The normal bundle of a substantial minimal surface $g\colon L^2\to\Q^N$ 
splits  orthogonally as the sum of plane bundles, except the last one which 
is a line bundle if the codimension is odd. 
These subbundles are spanned by the images of the higher fundamental forms.  
A minimal surface is called isotropic if the curvature ellipses of any order
but the highest one in odd codimension, that is, the images of the unit tangent circle 
under the higher fundamental forms, are circles at any point of $L^2$.  We observe
that isotropic surfaces have been studied under various names in quite 
different circumstances, namely, as superminimal, superconformal and 
pseudoholomorphic surfaces.

It is well known that any simply-connected minimal surface $g\colon L^2\to\Q^N$ 
in any codimension also allows a  one-parameter associated family  
$g_\theta\colon L^2\to\Q^N$, $\theta\in\Sf^1$, of isometric minimal immersions. 
The family is obtained integrating the system of Gauss, 
Codazzi and Ricci equations after rotating the second fundamental form a 
constant angle while keeping fixed the normal bundle and connection. 
If the surface is isotropic and lies in substantial
even codimension, then the associated family is trivial in the sense that the surfaces
are all pairwise congruent. 
On the other hand, in substantial odd codimension the associated family in 
never trivial.  Moreover, all surfaces in the associated family are isotropic 
and there are no other isotropic ones. 

Isotropic surfaces have been frequently considered  in the literature. 
The most obvious examples are holomorphic curves in $\C^N$ that can also be seen 
as isotropic surfaces in $\R^{2N}$. 
As for simply-connected isotropic surfaces in Euclidean 
space with odd codimension there is a Weierstrass type representation given 
in \cite{Cc} that generates all of them.  

When the ambient space is a round sphere, it results from the works of Calabi \cite{Cal} 
and Barbosa  \cite{B} that any topological $2$-sphere 
minimally immersed in $\Sf^N$ is isotropic and lies in even substantial codimension.
Miyaoka \cite{M} described the simply-connected isotropic surfaces in spheres 
with odd codimension in terms of the solutions of the affine Toda equation and, as
an application, managed to give all flat minimal tori in an explicit parametric form.  

It was shown by Vlachos \cite{V} that the so called Lawson's surfaces in spheres 
are isotropic but only for certain values of the parameters.
Examples of isotropic surfaces in $\Sf^5$ that include nonflat compact tori,
came out from Bryant's \cite{Br} theory  of pseudoholomorphic 
surfaces in the nearly Kaehler $\Sf^6$; see also \cite{BVW},\cite{BPW},\cite{EV} and \cite{H}.
Finally, observe that the polar surface (to be defined in \S $2$) associated to an 
isotropic substantial surface in an odd dimensional sphere is also isotropic.
\vspace{1ex}

Our first result is of local nature and follows from a quite simple argument.

\begin{theorem}\po\label{main1} Let $g\colon L^2\to\Q^{2n+1}$ 
be an oriented isotropic substantial surface. Then, the set of all  isotropic  
immersions of $L^2$ into $\Q^{2n+1}$ is either finite or a circle.
\end{theorem}

For the case of compact surfaces we have the main result of the paper.

\begin{theorem}\po\label{main2} Let $g\colon L^2\to\Sf^{2n+1}$
be an oriented isotropic substantial surface. If $L^2$ is compact,  then 
there exist at most finitely many  isotropic  immersions 
of $L^2$ into $\Sf^{2n+1}$.
\end{theorem}

The following is an immediate consequence of the last result.

\begin{corollary}\po\label{cor}
Let $g\colon L^2\to\Sf^{2n+1}$ be an oriented compact isotropic substantial 
surface. If $L^2$ admits a continuous one-parameter family of 
isometries $\varphi_t$ with $\varphi_0=id$, 
then there exists a continuous one-parameter family of isometries $\tau_t$ of  
$\Sf^{2n+1}$ such that $g\circ\varphi_t=\tau_t\circ g$ for
any value of the parameter.
\end{corollary}

\section{Preliminaries}

In this section, we collect several facts and definitions about minimal surfaces in space 
forms. For some of the details we refer to \cite{df1} and \cite{dv}. 
\vspace{1,5ex}

Let $g\colon L^2\to\Q^N$ denote an isometric immersion of a
two-dimensional  Riemannian manifold.
The $k^{th}$\emph{-normal space} of $g$ at $x\in L^2$ for $k\geq 1$ is 
defined as
$$
N^g_k(x)=\spa\{\a_g^{k+1}(X_1,\ldots,X_{k+1}):X_1,\ldots,X_{k+1}\in T_xL\}
$$
where 
$$
\a_g^s\colon TL\times\cdots\times TL\to N_gL,\;\; s\geq 3, 
$$
denotes the symmetric tensor called the $s^{th}$\emph{-fundamental form} given
inductively by
$$
\a_g^s(X_1,\ldots,X_s)=\left(\nabla^\perp_{X_s}\ldots
\nabla^\perp_{X_3}\a_g(X_2,X_1)\right)^\perp
$$
and $\a_g\colon TL\times TL\to N_gL$ stands for the standard second fundamental 
form of $g$ with values in the normal bundle.
Here  $\nabla^{\perp}$ denotes the induced connection in the normal bundle $N_gL$ of $g$ 
and $(\;\;)^\perp$ means taking the projection onto the normal complement of 
$N^g_1\oplus\ldots\oplus N^g_{s-2}$ in $N_gL$. 

Let $g\colon L^2\to\Q^N$ be a minimal isometric immersion. If $L^2$ is simply-connected, 
there exists a one-parameter \emph{associated family} of minimal isometric immersions 
$g_\theta\colon L^2\to\Q^N$, $\theta\in\Sf^1=[0,\pi)$ with real-analytic 
dependence on the parameter. To see this,  for each $\theta\in\Sf^1$ consider 
the orthogonal parallel tensor field 
$$
J_{\theta}=\cos\theta I+\sin\theta J
$$
where $I$ is the identity map.  Then, the symmetric section $\a_g(J_\theta\cdot,\cdot)$ 
of the bundle $\text{Hom}(TL\times TL,N_g L)$ satisfies the Gauss, Codazzi and Ricci 
equations with respect to the normal bundle and normal connection of $g$; 
see \cite{dg1} for details. 
Therefore, there exists an isometric minimal immersion  $g_{\theta}\colon L^2\to\Q^N$ 
whose second fundamental form is
\be\label{sff}
\a_{g_{\theta}}(X,Y)=\phi_\theta\a_g(J_{\theta}X,Y)
\ee 
where $\phi_\theta\colon N_gL\to N_{g_{\theta}}L$ is the parallel 
vector bundle isometry that identifies the normal subbundles
$N_j^g$ with $N_j^{g_\theta}$, $j\geq 1$.
\vspace{1ex}

A surface $g\colon L^2\to\Q^N$ is called \emph{regular} if for each $k$ the subspaces $N^g_k$ 
have constant dimension and thus form normal subbundles. 
Notice that regularity is always verified  along connected components of an open dense 
subset of $L^2$. 

Assume that an immersion $g\colon L^2\to\Q^N$ is minimal and substantial.  
The latter condition means that the codimension cannot be reduced.
In this case, the normal bundle of $g$ splits along the open dense subset 
of $L^2$ of regular points as
$$
N_gL=N_1^g\oplus N_2^g\oplus\dots\oplus N_m^g,\;\;\; m=[(N-1)/2],
$$
since all higher normal bundles  have rank two except possible the last 
one that has rank one if $N$ is odd; see \cite{Ch} or \cite{df1}. 
Moreover, if $L^2$ is oriented, then an orientation is induced on each plane 
vector bundle $N_s^g$  given by the ordered base
$$
\xi_1^s=\a_g^{s+1}(X,\ldots,X),\;\;\;\xi_2^s=\a_g^{s+1}(JX,\ldots,X)
$$
where $0\neq X\in TL$ and $J$ is the complex structure of $L^2$ 
determined by the metric and the orientation. 

If $g\colon L^2\to\Q^N$ is regular, then at $x\in L^2$ and for each 
$N_k^g$, $1\leq k\leq m$,
the \emph{$k^{th}$-order curvature ellipse}
$\E^g_k(x)\subset N^g_{k}(x)$ is defined  by
$$
\E^g_k(x) = \{\alpha_g^{k+1}(Z^{\varphi},\ldots,Z^{\varphi})\colon\,
Z^{\varphi}=\cos\varphi Z+\sin\varphi JZ\;\mbox{and}\;\varphi\in\Sf^1\}
$$
where $Z\in T_xL$ is any vector of unit length.

A regular  surface $g\colon L^2\to\Q^N$ is called \emph{isotropic} if it is  minimal 
and  at any $x\in L^2$ the ellipses of curvature $\E^g_k(x)$ contained in all 
two-dimensional  $N^g_k$$\,{}^{\prime}$s are circles.  
We point out that there are alternative ways to define 
isotropy for surfaces, for instance, in terms of the vanishing of higher order Hopf 
differentials \cite{V0} or as the mutual orthogonality of the complex line bundles 
determined by the harmonic sequences associated to the surface~\cite{BW}. 
\vspace{1,5ex}

\section{Isotropic surfaces}

In this section, we list several properties of isotropic surfaces in space forms 
already part of the literature. They may just justify the assumptions 
of our results or go further as to be ingredients of the proofs.
\vspace{1,5ex}

The first result follows from the arguments of Chern in \cite{Ch} and has been proved in
Proposition $4$ of \cite{V}.

\begin{fact}\po\label{fact1} If $g\colon L^2\to\Q^N$ is an  
isotropic surface, then the set $L_0$ where $g$ fails to be regular
is formed of isolated points and all $N^g_k$$\,{}^{\prime}$s extend smoothly to  $L_0$.
\end{fact}

 From either Corollary $6.1$ in \cite{J} or Corollary $5.2$ in \cite{V0} or Theorem $2$ 
in \cite{dv} we have the following result.

\begin{fact}\po\label{fact2} Any pair of isometric isotropic  surfaces  in $\Q^N$ that are
substantial in even codimension  are congruent.
\end{fact}

The situation for odd codimension is different since Theorem $2$ 
in \cite{dv} gives the following result for associated families.

\begin{fact}\po\label{fact3} Surfaces belonging to the associated family 
of a simply-connected isotropic surface in $\Q^N$ that is substantial in odd 
codimension are pairwise noncongruent.
\end{fact}

 For a simply-connected isotropic surface in odd substantial codimension 
first notice that it follows from the definition  that taking the associated family 
preserves isotropy. The following result is due to Johnson \cite{J}; 
see  Vlachos \cite{V0} for another proof.

\begin{fact}\po\label{fact4} If $g\colon L^2\to\Q^{2n+1}$ is an oriented 
simply-connected  isotropic substantial surface, then any other isotropic  substantial 
isometric immersion $f\colon L^2\to\Q^{2n+1}$ belongs to the associated 
family  $g_\theta$ of $g$. 
\end{fact}

Given a non-simply-connected  isotropic  oriented substantial surface 
$g\colon L^2\to\Q^{2n+1}$ let $f\colon L^2\to\Q^{2n+1}$ 
be another  isotropic  isometric substantial surface. 
Let $\pi\colon\tilde{L}^2\to L^2$ denote the universal covering map with 
$\tilde{L}^2$ equipped with the metric and orientation that makes $\pi$ an 
orientation preserving local isometry. In the sequel,  we denote
corresponding objects on $\tilde{L}^2$ with a tilde. 
Since the surfaces $\tilde{g}=g\circ \pi$ and $\tilde{f}=f\circ \pi$ 
are isotropic,  it follows from  Fact \ref{fact4} that  $\tilde{f}$ is congruent 
to some $\tilde{g}_{\theta}$ in the associated family of $\tilde{g}$. 
We thus have the following consequence.

\begin{fact}\po\label{fact5} Let $g\colon L^2\to\Q^{2n+1}$ be an oriented 
non-simply-connected isotropic substantial surface. Then, the set
$$
\mathcal{M}(g)=\left\{\theta\in\Sf^1\colon\text{there is}\;  
f\colon L^2\to\Q^{2n+1}\;\text{such that}\;f=\tilde{g}_{\theta}\circ\pi\right\}
$$ 
parametrizes the space of isometric isotropic substantial 
immersions of $L^2$ into $\Q^{2n+1}$.
\end{fact}

Let $g\colon L^2\to\Sf^{2n+1}$ be an isotropic substantial surface
and consider a smooth unit vector field  $e\in N_gL$ such that
$\spa\{e\}=N_n^g$. The surface 
$g^*=e\colon L^2\smallsetminus L_0\to\Sf^{2n+1}$ 
is usually called the \emph{polar} surface to $g$. The  first statement in
the next result follows from Proposition $8$ in \cite{df1} whereas the
second statement follows from Proposition $2$ and Theorem $2$ in \cite{V2}.

\begin{fact}\po\label{fact6} The polar surface 
$g^*\colon L^2\smallsetminus L_0\to\Sf^{2n+1}$  to a given oriented 
isotropic substantial surface $g\colon L^2\to\Sf^{2n+1}$ is also isotropic.
Moreover, the metric induced by $g^*$ is conformal to the 
metric of $L^2$ by a conformal factor that depends only on the metric of $L^2$.
\end{fact}

\section{The proofs}

In this section, we provide the proofs of our results 
stated in the introduction.
\vspace{1,5ex}

We first observe  the following elementary fact.

\begin{proposition}\label{D}\po Let $g\colon L^2\to\Q^N$
be a non-simply-connected  isotropic surface. 
For any $\sigma\in\mathcal{D}$ in the group of deck transformations of 
$\pi\colon\tilde{L}^2\to L^2$ the surfaces 
$\tilde{g}_{\theta}\colon\tilde{L}^2\to\Q^N$ and 
$\tilde{g}_{\theta}\circ\sigma\colon\tilde{L}^2\to\Q^N$ are congruent
for any $\theta\in \Sf^1$.
\end{proposition}

\proof It is sufficient to show the existence of a parallel vector bundle 
isometry between the normal bundles of 
$\tilde{g}_{\theta}$ and $\tilde{g}_{\theta}\circ \sigma$ that preserves 
the second fundamental forms. Let $\phi_{\theta}$ be the vector isometry
between  the normal bundles of $\tilde{g}$ and $\tilde{g}_{\theta}$. Define 
a vector bundle isometry 
${\psi}_{\theta}\colon N_{\tilde{g}_{\theta}}\tilde{L} \to
N_{\tilde{g}_{\theta}\circ \sigma}\tilde{L}$ by
$$
{\psi}_{\theta}\xi= \phi_{\theta}(\eta \circ \sigma^{-1})\circ \sigma
$$
where  $\xi=\phi_{\theta}\eta$ for $\eta\in N_{\tilde{g}}\tilde L$. 
The second fundamental forms of $\tilde{g}_{\theta}\circ \sigma$ and 
$\tilde{g}_{\theta}$ relate by
$$
\a_{\tilde{g}_{\theta}\circ\sigma}(\tilde{X}, \tilde{Y})
=\phi_{\theta}\a_{\tilde{g}}\big(\tilde{J}_{\theta} 
\sigma_*\tilde{X},\sigma_*\tilde{Y}\big)  
$$
for any $\tilde{X},\tilde{Y}\in T\tilde{L}$, where 
$$
\tilde{J}_{\theta}=\cos\theta\tilde{I}+\sin\theta\tilde{J}.
$$
Since $\sigma$ is a deck transformation, then
$$
\a_{\tilde{g}_{\theta}\circ \sigma}
={\psi}_{\theta}\circ\a_{\tilde{g}_{\theta}}.
$$
Let $\xi=\phi_{\theta}\eta$ where $\eta\in N_{\tilde{g}}\tilde L$. 
We have that
\bea
(\nabla^{\perp}_{\tilde{X}}\mathit{\psi}_{\theta})\xi
\!\!\!&=&\!\!\!
\nabla^{\perp}_{\tilde{X}}\big(\phi_{\theta}(\eta\circ\sigma^{-1})\circ\sigma  
\big)-\phi_{\theta}\big(\nabla^{\perp}_{\tilde{X}}(\eta\circ\sigma^{-1})\big)\circ\sigma\\
\!\!\!&=&\!\!\!
\big(\nabla^{\perp}_{\sigma_*\tilde{X}} \phi_{\theta}(\eta\circ\sigma^{-1})\big)\circ\sigma  
-\phi_{\theta}\big(\nabla^{\perp}_{\tilde{X}}(\eta\circ \sigma^{-1})\big)\circ\sigma\\
\!\!\!&=&\!\!\!
\phi_{\theta}\big( \nabla^{\perp}_{\sigma_*\tilde{X}}(\eta\circ\sigma^{-1})  
-\nabla^{\perp}_{\tilde{X}}(\eta\circ\sigma^{-1})\big)\circ\sigma
\eea
where $\tilde{X}\in T\tilde{L}$ and $\nabla^{\perp}$ stands for the normal 
connection of $\tilde{g}_{\theta}$ as well as for $\tilde{g}_{\theta}\circ\sigma$. 
Now let  $\delta\in N_g L$ be such that 
$\eta\circ\sigma^{-1}=\delta\circ\pi$.  Observe that
$$
\nabla^{\perp}_{\sigma_*\tilde{X}}\eta\circ\sigma^{-1}  
-\nabla^{\perp}_{\tilde{X}}\eta\circ\sigma^{-1}
=\nabla^{\perp}_{\pi_*\sigma_*\tilde{X}} \delta  
-\nabla^{\perp}_{\pi_*\tilde{X}}\delta =0,
$$
and thus $\psi_\theta$ is parallel.\vspace{1,5ex}\qed

\noindent\emph{Proof of Theorem \ref{main1}:} Proposition \ref{D} yields 
a map $\varPhi_{\theta}\colon\mathcal{D}\to\text{Isom}(\Q^{2n+1})$ 
for each given $\theta\in\Sf^1$ such that
$$
\tilde{g}_{\theta}\circ\sigma=\varPhi_{\theta}(\sigma)\circ\tilde{g}_{\theta}.
$$
Thus $\theta\in\mathcal{M}(g)$ if and only if 
$\varPhi_{\theta}(\mathcal{D})=\{I\}$.
Suppose that $\mathcal{M}(g)$ is not finite and let $\{\theta_m\}_{m\in\mathbb{N}}$ 
be a sequence in $\mathcal{M}(g)$  that  converges to some $\theta_0\in\Sf^1$. 
Hence $\varPhi_{\theta_0}(\mathcal{D})=\{I\}$. 

Take $\sigma\in\mathcal{D}$.
By the Mean value Theorem applied to each entry
$(\varPhi_{\theta}(\sigma))_{jk}$ of the corresponding matrix, we have
\be\label{nul}
\dfrac{d}{d\theta}(\varPhi_{\theta}(\sigma))_{jk}(\mathring{\theta}_m)=0
\ee
for some $\mathring{\theta}_m$ which lies between $\theta_0$ and $\theta_m$.
By continuity it follows that
$$
\dfrac{d}{d\theta}(\varPhi_{\theta}(\sigma))_{jk}(\theta_0)=0.
$$
Consider the sequence $\{\mathring{\theta}_m\}_{m\in\mathbb{N}}$ 
that converges $\theta_0$ and observe that in view of (\ref{nul}) similar
argument gives
$$
\dfrac{d^2}{d\theta^2}(\varPhi_{\theta}(\sigma))_{jk}(\theta_0)=0.
$$
Repeating the argument yields 
$$
\dfrac{d^n}{d\theta^n}(\varPhi_{\theta}(\sigma))_{jk}(\theta_0)=0
$$
for any integer $n\geq 1$. Since $\varPhi_{\theta}(\sigma)$ is an
analytic curve  in $\text{Isom}(\Q^{2n+1})$, we conclude that 
$\varPhi_{\theta}(\sigma)=I$ for each 
$\sigma\in\mathcal{D}$ and Fact \ref{fact5} now gives
$\mathcal{M}(g)=\Sf^1$.\qed\vspace{1,5ex}

The proof of the following fact is immediate.

\begin{lemma}\label{isomor}\po
Let $g\colon L^2\rightarrow\Q^{2n+1}$ be an isotropic substantial surface. 
For any $\theta \in \mathcal{M}(g)$ there exists a parallel vector bundle isometry 
$\phi_{\theta}\colon N_g L\to N_{\tilde{g}_{\theta}\circ\pi}L$ such that the higher fundamental 
forms of $g$  and $\tilde{g}_{\theta}\circ\pi$ are related by
$$
\a^s_{\tilde{g}_{\theta}\circ\pi}(X_1,\dots,X_s)=\phi_{\theta}\a^s_g(J_{\theta}X_1,\dots,X_s)
$$
where $X_1,\dots,X_s\in TL$.
\end{lemma}

\noindent\emph{Proof of Theorem \ref{main2}:}
Choose a positively oriented orthonormal frame $\{e_1,e_2\}$ in $TL$ along
an open subset of $L^2\smallsetminus L_0$. Since $g$ is isotropic, we may choose 
a local orthonormal normal frame $\{e_3,\ldots,e_{2n+1}\}$ such that
$\{e_{2s+1},e_{2s+2}\}$ is given by
 $$
\a^{s+1}_g(e_1,\dots,e_1, e_1)=\kappa_s e_{2s+1}\;\;\mbox{and}\;\;
\a^{s+1}_g(e_1,\dots,e_1, e_2)=\kappa_s e_{2s+2},\;\; 1\leq s \leq n-1
$$
where $\kappa_s$ denotes the radii  of the circular ellipse of curvature $\E^g_s$.
Then $e_{2n+1}$ spans the last normal bundle.  
According to Lemma 5 in \cite{V} the normal connection forms 
$$
\omega_{\a\beta}=\< \nap e_{\a},e_{\beta}\>
$$
satisfy
\be\label{con1}
\omega_{2s,2s+1}=-*\omega_{2s-1,2s+1},\;\;\omega_{2s,2s+2}=-*\omega_{2s-1,2s+2},
\ee
\be\label{con2}
\omega_{2s-1,2s+2}=*\omega_{2s-1,2s+1},\;\;\omega_{2s,2s+2}=*\omega_{2s,2s+1},
\;\; 1\leq s \leq n-1,
\ee
and
\be\label{con3}
\omega_{2n,2n+1}=-*\omega_{2n-1,2n+1}
\ee
where $*$ denotes the Hodge operator, i.e., $*\omega=-\omega\circ J$.

The metric and second fundamental of $g$ can be complex linearly extended to the 
complexified tangent bundle $TL\otimes\mathbb{C}$ and complexified normal 
bundle $N_gL\otimes\mathbb{C}$.
Setting $E=e_1-ie_2$, then (\ref{con1}), (\ref{con2}) and  (\ref{con3}) yield

\be\label{coni}
\omega_{2s-1,2s+2}(E)=-i\omega_{2s-1,2s+1}(E),
\;\;\omega_{2s,2s+1}(E)=i\omega_{2s-1,2s+1}(E),
\ee
\be\label{conii}
\;\omega_{2s,2s+2}(E)=\omega_{2s-1,2s+1}(E),\;\; 1\leq s \leq n-1,
\ee
and 
\be\label{coniii}
\;\omega_{2n,2n+1}(E)=i\omega_{2n-1,2n+1}(E).
\ee

Take $\theta_1<\dots<\theta_m\in\mathcal{M}(g)$. In the sequel, we 
regard the surfaces $g_{\theta_j}=\tilde{g}_{\theta_j}\circ\pi$ as 
lying in $\R^{2n+2}$ by composing with the 
umbilical inclusion $\Sf^{2n+1}\hookrightarrow\R^{2n+2}$. 
We claim that if
\be\label{coor}
\sum_{j=1}^m\<g_{\theta_j},v_j\>=0
\ee
for some $v_j\in\R^{2n+2}$, then $v_j=0$ for all $1\leq j\leq m$.

To prove the claim, we may assume that $v_j\neq 0$ for all $1\leq j\leq m$. 
Differentiating (\ref{coor}) yields
\be\label{d}
\sum_{j=1}^m\<g_{{\theta_j}_*}, v_j\>=0 \;\; \mbox{and}\;\; 
\sum_{j=1}^m\<\a_{g_{\theta_j}}, v_j\> =0.
\ee
Then, Lemma \ref{isomor} and the second equation in (\ref{d}) give
$$
\sum_{j=1}^m\<\phi_{\theta_j}\a_g(J_\theta E,E),v_j\>=0.
$$
Since $J_\theta E=e^{i\theta} E$, it follows easily that
\be\label{a}
\sum_{j=1}^m e^{i\theta_j}\<\phi_{\theta_j}(e_3-ie_4),v_j\>=0.
\ee
Differentiating  (\ref{a}) with respect to $X \in TM\otimes \mathbb{C}$ yields
$$
\sum_{j=1}^m e^{i\theta_j}\<g_{{\theta_j}_*} A_{\phi_{\theta_j}(e_3-ie_4)}X,v_j\>
-\sum_{j=1}^m e^{i\theta_j}\<\nap_X\phi_{\theta_j}(e_3-ie_4),v_j\>=0
$$
where  $A_{\phi_{\theta_j}\eta}$ is the shape operator of $g_{\theta_j}$
in the direction $\phi_{\theta_j}\eta$. Lemma \ref{isomor} and (\ref{sff}) give
$$
A_{\phi_{\theta_j}(e_3-ie_4)} =e^{-i\theta_j}A_{e_3-ie_4}
$$
where $A_\eta$ denotes the shape operator of $g$. This and the 
first equation in (\ref{d}) yield
$$
\sum_{j=1}^m e^{i\theta_j}\<\nap_X\phi_{\theta_j}(e_3-ie_4),v_j\>=0.
$$
In view of  (\ref{a}) this equation is equivalent to
\be\label{56}
(\omega_{35}(X)-i\omega_{45}(X))\sum_{j=1}^m e^{i\theta_j}\<\phi_{\theta_j}e_5,v_j\> +
(\omega_{36}(X)-i\omega_{46}(X))\sum_{j=1}^m e^{i\theta_j}\<\phi_{\theta_j}e_6,v_j\>=0.
\ee
On the other hand,  we obtain from (\ref{coni}) and (\ref{conii}) that
\be\label{35}
\omega_{45}(E)=i\omega_{35}(E),\;\;\omega_{36}(E)=-i\omega_{35}(E) 
\;\;\mbox{and}\;\;\omega_{46}(E)=\omega_{35}(E).
\ee
Using (\ref{35}) it follows that (\ref{56}) for $X=E$ becomes
$$
\omega_{35}(E)\sum_{j=1}^me^{i\theta_j}\<\phi_{\theta_j}(e_5-ie_6),v_j\>=0.
$$
We have from (\ref{35}) that $\omega_{35}(E)$ cannot vanish 
since $g$ is substantial. Hence, 
$$
\sum_{j=1}^me^{i\theta_j}\<\phi_{\theta_j}(e_5-ie_6),v_j\>=0.
$$

We argue that 
\be\label{r}
\sum_{j=1}^m e^{i\theta_j}\<\phi_{\theta_j}(e_{2s+1}-ie_{2s+2}),v_j\>=0,\;\;s\leq n-1,
\ee
and
\be\label{rr}
\sum_{j=1}^me^{i\theta_j}\<\phi_{\theta_j}e_{2n+1},v_j\>=0.
\ee
By induction, we  assume (\ref{r}) for $s=r$  
and differentiate  with respect to $E$. We obtain
\bea
(i\omega_{2r-1,2r+2}(E)\!\!\!&-&\!\!\!\omega_{2r-1,2r+1}(E))
\sum_{j=1}^m e^{i\theta_j}\<\phi_{\theta_j}e_{2r-1},v_j\>\\
\!\!\!&+&\!\!\!(i\omega_{2r,2r+2}(E)-\omega_{2r,2r+1}(E))
\sum_{j=1}^m e^{i\theta_j}\<\phi_{\theta_j}e_{2r},v_j\>\\
\!\!\!&+&\!\!\!i\omega_{2r+1,2r+2}(E)
\sum_{j=1}^m e^{i\theta_j}\<\phi_{\theta_j}(e_{2r+1}-ie_{2r+2}),v_j\>\\
\!\!\!&+&\!\!\!(\omega_{2r+1,2r+3}(E)-i\omega_{2r+2,2r+3}(E))
\sum_{j=1}^m e^{i\theta_j}\<\phi_{\theta_j}e_{2r+3},v_j\>\\
\!\!\!&+&\!\!\!(\omega_{2r+1,2r+4}(E)-i\omega_{2r+2,2r+4}(E))
\sum_{j=1}^m e^{i\theta_j}\<\phi_{\theta_j}e_{2r+4},v_j\>\\
\!\!\!&=&\!\!\!0
\eea
if $r\leq n-2$.  If $r=n-1$, we have
\bea
(i\omega_{2n-3,2n}(E)\!\!\!&-&\!\!\!\omega_{2n-3,2n-1}(E))
\sum_{j=1}^m e^{i\theta_j}\<\phi_{\theta_j}e_{2n-3},v_j\>\\
\!\!\!&+&\!\!\!(i\omega_{2n-2,2n}(E)-\omega_{2n-2,2n-1}(E))
\sum_{j=1}^me^{i\theta_j}\<\phi_{\theta_j}e_{2n-2},v_j\>\\
\!\!\!&+&\!\!\!i\omega_{2n-1,2n}(E)
\sum_{j=1}^m e^{i\theta_j}\<\phi_{\theta_j}(e_{2n-1}-ie_{2n}),v_j\>\\
\!\!\!&+&\!\!\!(\omega_{2n-1,2n+1}(E)-i\omega_{2n,2n+1}(E))
\sum_{j=1}^m e^{i\theta_j}\<\phi_{\theta_j}e_{2n+1},v_j\>\\
\!\!\!&=&\!\!\!0.
\eea
Using (\ref{coni}), (\ref{conii}) and (\ref{coniii}), we obtain
$$
\omega_{2r+1,2r+3}(E)\sum_{j=1}^m e^{i\theta_j}
\<\phi_{\theta_j}(e_{2r+3}-ie_{2r+4}),v_j\>=0,\;\;r\leq n-2,
$$
and 
$$
\omega_{2n-1,2n+1}(E)\sum_{j=1}^m e^{i\theta_j}
\<\phi_{\theta_j}e_{2n+1},v_j\>=0.
$$
Since $g$ is substantial, from (\ref{coni}), (\ref{conii}) and (\ref{coniii}) 
we see that $\omega_{2r+1,2r+3}(E)$, $r\leq n-1$, cannot vanish. This completes 
the inductive argument and proves (\ref{r}) and (\ref{rr}).

>From Fact \ref{fact6} the polar surface $g^*_\theta$ to $g_\theta$ 
is isometric to $g^*$ for any $\theta\in\mathcal{M}(g)$.
It follows from (\ref{rr}) that the polar surfaces 
$g^*_{\theta_j}=\phi_{\theta_j}e_{2n+1}\colon L^2\smallsetminus L_0\to\Sf^{2n+1}$
satisfy
$$
\sum_{j=1}^m e^{i\theta_j}\<g^*_{\theta_j},v_j\>=0.
$$
We easily see that
\be\label{pm}
\sum_{j=1}^{m-1}\<g^*_{\theta_j},v_j^1\>=0
\ee
where $v_j^1=\sin(\theta_m-\theta_j)v_j\neq 0$.
We have from Fact \ref{fact6} that the surfaces $g^*_{\theta_j}$, 
$1\leq j\leq m-1$, are isotropic. 
If $m>2$ and since the polar surface of $g^*_{\theta_j}$ is just 
$g_{\theta_j}$, arguing  as for (\ref{pm}) we obtain that
$$
\sum_{j=1}^{m-2}\<g_{\theta_j},v_j^2\>=0
$$
where $v_j^2=\sin(\theta_{m-1}-\theta_j)v_j^1\neq 0$.
By repeating the argument, if necessary, we conclude that either 
$\<g_{\theta_1},v\>=0$ or $\<g^*_{\theta_1},v\>=0$  for some $v\neq 0$. 
But from the definition of the associated family we have that any 
$g_\theta$ is substantial. That any $g^*_{\theta}$ must also be substantial follows 
easily using Proposition $8$ in \cite{df1}.  Thus, we reached a contradiction
and the claim has  been proved.

We now conclude the proof. 
According to Theorem \ref{main1} the set $\mathcal{M}(g)$ is either finite 
or $\Sf^1$. Suppose that $\mathcal{M}(g)=\Sf^1$.  Then any finite subset of
height functions of the surfaces $g_{\theta},\,\theta\in\Sf^1$, must be linearly 
independent. But since these functions are eigenfunctions of the Laplace
operator on  $L^2$ with eigenvalue $2$, this  contradicts the fact that any
eigenspace of the Laplace operator on a compact manifold  has finite dimension.
\vspace{1,5ex}\qed

\noindent\emph{Proof of Corollary \ref{cor}:}
 From Theorem \ref{main2} we know that $\mathcal{M}(g)$ is finite.
Consider the isotropic immersions  $g^t=g\circ\varphi_t$.  
Since the second fundamental form of $g^t$ depends continuously on the 
parameter, we deduce that $g^t$ for any $t$ is congruent to exactly one 
$g_\theta$ for some $\theta\in\mathcal{M}(g)$. Since $\varphi_0=id$, 
by continuity we conclude that $g^t$ is congruent to $g$ for any $t$.\qed

%\vspace{.5in} 
{\renewcommand{\baselinestretch}{1}
\hspace*{-20ex}\begin{tabbing} \indent\= IMPA -- Estrada Dona Castorina, 110
\indent\indent\= Univ. of Ioannina -- Math. Dept. \\
\> 22460-320 -- Rio de Janeiro -- Brazil  \>
45110 Ioannina -- Greece \\
\> E-mail: marcos@impa.br \> E-mail: tvlachos@uoi.gr
\end{tabbing}}
\end{document}